\documentclass[12pt,english,twoside]{article}
\usepackage[english]{babel}
\usepackage{amsmath,amsfonts,amssymb}
\usepackage{mathtools}
\usepackage{color}
\usepackage{natbib}
\setcitestyle{numbers,square}
\usepackage[cp850]{inputenc}
\usepackage{tikz}

\allowdisplaybreaks
\newcounter{theorem}[section]
\renewcommand{\thetheorem}{\thesection.\arabic{theorem}}
\newenvironment{lemma}[      1]{\refstepcounter{theorem} %
\bf \thetheorem\ Lemma#1.        \it}{}
\newenvironment{theorem}[    1]{\refstepcounter{theorem} %
\bf \thetheorem\ Theorem#1.      \it}{}
\newenvironment{proposition}[1]{\refstepcounter{theorem} %
\bf \thetheorem\ Proposition#1.  \it}{}
\newenvironment{corollary}[  1]{\refstepcounter{theorem} %
\bf \thetheorem\ Corollary#1.    \it}{}
\newenvironment{example}[    1]{\refstepcounter{theorem} %
\bf \thetheorem\ Example#1.      \rm}{\hfill $ \Box $}
\newenvironment{examples}[   1]{\refstepcounter{theorem} %
\bf \thetheorem\ Examples#1.     \rm}{\hfill $ \Box $}
\newenvironment{remark}[  1]{\refstepcounter{theorem} %
\bf \thetheorem\ Remark#1.       \rm}{}

\newenvironment{proof}{ %
\it              Proof.          \rm}{\hfill $ \Box $}

\newcounter{abc}[theorem]

\newcounter{one}[theorem]
\newenvironment{onelist}{\begin{list}{%
\rm (\arabic{one}) \hfill           }{\usecounter{one} %
\topsep0mm \partopsep0mm \parsep0mm \itemsep0mm %
\leftmargin2em \labelwidth2em \labelsep0em}}{\end{list}}
\newenvironment{1onelist}{\begin{list}{%
\rm (1.\arabic{one}) \hfill           }{\usecounter{one} %
\topsep0mm \partopsep0mm \parsep0mm \itemsep0mm %
\leftmargin2em \labelwidth2em \labelsep0em}}{\end{list}}
\newenvironment{2onelist}{\begin{list}{%
\rm (2.\arabic{one}) \hfill           }{\usecounter{one} %
\topsep0mm \partopsep0mm \parsep0mm \itemsep0mm %
\leftmargin2em \labelwidth2em \labelsep0em}}{\end{list}}

\newenvironment{hylist}{\begin{list}{%
--               \hfill            }{%
\topsep0mm \partopsep0mm \parsep0mm \itemsep0mm %
\leftmargin2em \labelwidth2em \labelsep0em}}{\end{list}}

\newcounter{rom}
\newenvironment{romlist}{\begin{list}{%
\rm (\roman{rom})  \hfill           }{\usecounter{rom} %
\topsep0mm \partopsep0mm \parsep0mm \itemsep0mm %
\leftmargin2em \labelwidth2em \labelsep0em}}{\end{list}}

\newcommand{\I}{{\mathbb I}}
\newcommand{\N}{{\mathbb N}}
\newcommand{\R}{{\mathbb R}}
\newcommand{\Real}{{\mathbb R}}

\newcommand{\BB}{{\cal B}}
\newcommand{\CC}{{\cal C}}
\newcommand{\DD}{{\cal D}}
\newcommand{\EE}{{\cal E}}

\newcommand{\aaa}{{\bf a}}
\newcommand{\bbb}{{\bf b}}

\newcommand{\eee}{{\bf e}}
\newcommand{\uuu}{{\bf u}}
\newcommand{\vvv}{{\bf v}}

\newcommand{\UUU}{{\bf U}}

\newcommand{\eeta}{{\boldsymbol\eta}}

\newcommand{\zero}{{\bf 0}}
\newcommand{\eins}{{\bf 1}}

\begin{document}
\title{\Large\bf On Minimal Copulas under the Concordance Order}
\author{\normalsize Jae Youn Ahn\footnote{E--mail Jae Youn Ahn: {\tt jaeyahn@ewha.ac.kr} (Corresponding author)} \\ [1ex]
\small Ewha Womans University  \\[-0.8ex]
\small Seoul, Korea
\and 
\normalsize Sebastian Fuchs\footnote{E--mail Sebastian Fuchs: {\tt sfuchs@statistik.tu-dortmund.de} (Corresponding author)} \\[1ex]
\small Technische Universit{\"a}t Dortmund \\[-0.8ex] 
\small Dortmund, Germany }
\date{}
\maketitle

\bigskip

\begin{abstract}
\noindent 
In the present paper, we study extreme negative dependence focussing on the concordance order for copulas.
With the absence of a least element for dimensions $d\geq3$, the set of all minimal elements in the collection of all copulas turns out to be a natural and quite important extreme negative dependence concept.
We investigate several sufficient conditions and we provide a necessary condition for a copula to be minimal:
The sufficient conditions are related to the extreme negative dependence concept of $d-$countermonotonicity 
and the necessary condition is related to the collection of all copulas minimizing multivariate Kendall's tau.
The concept of minimal copulas has already been proved to be useful in various continuous and concordance order preserving optimization problems including variance minimization and the detection of lower bounds for certain measures of concordance.
We substantiate this key role of minimal copulas by showing that every continuous and concordance order preserving functional on copulas is minimized by some minimal copula and, in the case the continuous functional is even strictly concordance order preserving, it is minimized by minimal copulas only.
Applying the above results, we may conclude that every minimizer of Spearman's rho is also a minimizer of Kendall's tau.
\end{abstract}




\section{Introduction}

\smallskip

The strongest notion of positive dependence is given by comonotonicity. 
A $d-$dimensional continuous random vector is said to be \emph{comonotonic} if one of the following equivalent conditions is satisfied; see, e.g., \cite{dus2016}:
\begin{hylist}
\item
the random vector has the upper Fr{\'e}chet--Hoeffding bound $M$ as its copula,

\item 
each of its coordinates is almost surely an increasing transformation of the other,
 
\item
each of its bivariate subvectors is comonotonic. 
\end{hylist}
In many applications in finance and insurance, comonotonicity can be regarded as one of the most dangerous behaviours and, \pagebreak
with its popularity in risk management as an extreme positive dependence, 
it seems quite reasonable that also extreme negative dependence is getting more attention, see, e.g., \citep{Dhaene}.
In the bivariate case, 
the strongest notion of negative dependence is given by countermonotonicity. 
A bivariate continuous random vector is said to be \emph{countermonotonic} if one of the following equivalent conditions is satisfied; see, e.g., \cite{dus2016}:
\begin{hylist}
\item
the random vector has the lower Fr{\'e}chet--Hoeffding bound $W$ as its copula,

\item 
each of its coordinates is almost surely a decreasing transformation of the other.
\end{hylist}
In contrast to comonotonicity, for dimensions $d\geq 3$ there exist no continuous random vector for which all the bivariate subvectors are countermonotonic and hence there exists no single agreed definition of extreme negative dependence in arbitrary dimension;
see, e.g., \cite{Cheung4, Dhaene4, fumcsc2018, Ahn7, wawa2011, Ruodu4}.



\bigskip

In the present paper, we investigate extreme negative dependence focussing on the concordance order for copulas.
It is well--known that
$M$ is the greatest element
and, in the bivariate case, 
$W$ is the least element
in the collection of all copulas with respect to concordance order.
With the absence of a least element for dimensions $d\geq3$, the set of all minimal elements
(so--called \emph{minimal copulas}), 
i.e. all locally least elements, in the collection of all copulas turns out to be a natural and quite important extreme negative dependence concept, 
in particular, with regard to the minimization of continuous and concordance order preserving optimization problems.
This includes several measures of dependence like Kendall's tau and Spearman's rho,
but also the variance of the sum of several given random variables as a map on copulas;
see, e.g., \cite{daniel14}.
While it is well--known that every continuous and (strictly) concordance order preserving functional on copulas 
is (uniquely) maximized by $M$ and, 
in the bivariate case, (uniquely) minimized by $W$,
less is known about minimization for dimensions $d \geq 3$. 

Most recently, \citet{Ahn11} and \citet{Ahn9} have demonstrated the potential of minimal copulas in variance minimization when the marginals are uniform, elliptical or belong to the unimodal--symmetric location--scale family;
we refer to \cite{Ruschendorf2, wawa2011} for further results on variance minimization.
Moreover, \citet{genebg2011} and \citet{Ahn7} have provided minimal copulas minimizing certain measures of concordance.
Further continuous and concordance order preserving functionals on copulas are discussed in 
\cite{daniel14, fuc2016bic}, but also in \cite{luxpa2017, prei2016, ruodu3, Tankov} where the necessary properties of the functional strictly depend on the function to be integrated; 
for more details on this topic we refer to \cite{ruodu3} and the references therein.

\bigskip

In this paper, we first discuss the existence and list some important examples of copulas that are minimal with respect to concordance order. 
We then investigate several sufficient conditions (Proposition \ref{dCMMinimal} and Theorem \ref{dCMProduct}) and we provide a necessary condition (Theorem \ref{MinimalKendallCM}) for a copula to be minimal. 
The sufficient conditions are related to the extreme negative dependence concept of $K-$countermonotonicity introduced in \cite{Ahn7,Ahn9}
and the necessary condition is related to the collection of all copulas minimizing multivariate Kendall's tau:
It turns out that every minimal copula minimizes Kendall's tau which is the main result of this paper.
We further point out the key role of minimal copulas with regard to the minimization of continuous and (strictly) concordance order preserving optimization problems: 
It turns out that every continuous and concordance order preserving functional on copulas is minimized by some minimal copula and that, in the case the continuous functional is even strictly concordance order preserving, it is minimized by minimal copulas only (Theorem \ref{FunctionalCOP.1}). 
Finally, we apply our results to Kendall's tau and Spearman's rho and show that every minimizer of Spearman's rho is also a minimizer of Kendall's tau (Corollary \ref{MOCKendallSpearman}).


\section{Preliminaries}

\smallskip

In this section, 
we fix some notation and recall some definitions and results on 
copulas, a group of transformations of copulas and the concordance order.

\bigskip

Let   $ \I := [0,1] $  
and let   $ d\geq2 $   be an integer which will be kept fix throughout this paper. 
We denote 
by   $ \eee_1,\dots,\eee_d $   the unit vectors in   $ \R^d $, 
by   $ \zero $   the vector in   $ \R^d $   with all coordinates being equal to   $ 0 $   and 
by   $ \eins $   the vector in   $ \R^d $   with all coordinates being equal to   $ 1 $.
For   $ \uuu,\vvv\in\R^d $, 
we use the notation   $ \uuu\leq\vvv $   resp.\   $ \uuu<\vvv $   in the usual sense 
such that   $ u_k \leq v_k $   resp.\   $ u_k<v_k $   holds for every   $ k\in\{1,\dots,d\} $.

\subsection*{Copulas}

For   $ K\subseteq\{1,...,d\} $, 
we consider the map   $ \eeta_K : \I^d\times\I^d \to \I^d $   given coordinatewise by
$$
  (\eeta_K(\uuu,\vvv))_k                
 := \begin{cases}                         
         u_k, & k \in \{1,...,d\} \setminus K,  \\
         v_k, & k \in K,                        
        \end{cases}                           
$$
and for $ k \in \{1,...,d\} $   we put   $ \eeta_k := \eeta_{\{k\}} $. 
A 
\emph{copula} is a function   $ C: \I^{d} \to \I $   satisfying the following conditions:
\begin{romlist}
\item
$ \sum_{K \subseteq \{1,...,d\}} (-1)^{d-|K|}\,C(\eeta_K(\uuu,\vvv)) \geq 0  $   
holds for all   $ \uuu,\vvv\in\I^d $   such that   $ \uuu\leq\vvv $. 

\item
$ C(\eeta_k(\uuu,\zero)) = 0   $   holds for every   $ k\in\{1,...,d\} $   and every   $ \uuu\in\I^d $. 

\item  
$ C(\eeta_k(\eins,\uuu))	= u_k $   holds for every   $ k\in\{1,...,d\} $   and every   $ \uuu\in\I^d $. 
\end{romlist}  
This definition is in accordance with the literature;
see, 
e.\,g., 
\cite{dus2016, nel2006}.
The collection   $ \CC $   of all copulas is convex. 
The following copulas are of particular interest:

\begin{hylist}
\item   
The \emph{upper Fr{\'e}chet--Hoeffding bound} $M$ given by 
$ M(\uuu) := \min\{u_1,\dots,u_d\} $ 
is a copula and every copula $C$ satisfies $C(\uuu) \leq M(\uuu)$ for every $\uuu \in \I^d$. 

\item   
The \emph{product copula}   $ \Pi $   given by 
$ \Pi(\uuu) := \prod_{k=1}^d u_k $ 
is a copula. 

\item   
The \emph{lower Fr{\'e}chet--Hoeffding bound} $W$ given by 
$ W(\uuu) := \max\{ \sum_{k=1}^d u_k + 1 - d, 0 \} $   
is a copula only for $d=2$, and every copula $C$ satisfies $W(\uuu) \leq C(\uuu)$ for every $\uuu \in \I^d$.  
\end{hylist}
Since every copula   $ C $   has a unique extension to a distribution function   $ \R^d\to\I $, 
there exists a unique probability measure   $ Q^C : \BB(\I^d)\to\I $   satisfying   
$ Q^C[[\zero,\uuu]] = C(\uuu) $
for every   $ \uuu\in\I^d $. 
The probability measure   $ Q^C $   is said to be the 
\emph{copula measure} with respect to   $ C $   and it satisfies   $ Q^C[(\uuu,\vvv)] = Q^C[[\uuu,\vvv]] $   for all   $ \uuu,\vvv\in\I^d $   
such that   $ \uuu\leq\vvv $.

\subsection*{A Group of Transformations of Copulas}

Let   $ \Phi $   denote the collection of all transformations   $ \CC\to\CC $   
and consider the composition   $ \circ: \Phi\times\Phi \to \Phi $   given by   $ (\varphi_2\circ\varphi_1)(C) := \varphi_2(\varphi_1(C)) $   
and the map   $ \iota\in\Phi $   given by   $ \iota(C) := C $. 
Then   $ (\Phi,\circ) $   is a semigroup with neutral element   $ \iota $. 
For   $ i,j,k\in\{1,...,d\} $   with   $ i \neq j $, 
we define the maps   $ \pi_{i,j},\nu_k : \CC\to\CC $   by letting 
\begin{eqnarray*}
        (\pi_{i,j}(C))(\uuu)                                      
& := &  C(\eeta_{\{i,j\}}(\uuu,u_j\,\eee_i+u_i\,\eee_j))          \\*
        (\nu_k    (C))(\uuu)                                      
& := &  C(\eeta_k(\uuu,\eins)) - C(\eeta_k(\uuu,\eins\!-\!\uuu))  
\end{eqnarray*}
Each of these maps is an involution and there exists 
\begin{hylist}
\item   a smallest subgroup   $ \Gamma^\pi $   of   $ \Phi $   containing every   $ \pi_{i,j} $, 
\item   a smallest subgroup   $ \Gamma^\nu $   of   $ \Phi $   containing every   $ \nu_k $   and 
\item   a smallest subgroup   $ \Gamma     $   of   $ \Phi $   containing   $ \Gamma^\pi\cup\Gamma^\nu $. 
\end{hylist}
The transformations in $\Gamma^{\pi}$ are called \emph{permutations} and the transformations in $\Gamma^\nu$ are called \emph{reflections}.
The group   $ \Gamma^\nu $   is commutative and for   $ K \subseteq \{1,...,d\} $   we define 
$$ 
  \nu_K := \bigcirc_{k \in K} \nu_{k} 
$$
(such that   $ \nu_\emptyset = \iota $). 
We note that the 
\emph{total reflection} 
$\tau := \nu_{\{1,...,d\}}$
transforms every copula into its \emph{survival copula}.
From a probabilistic viewpoint, if $\UUU$ is a random vector whose distribution function is the copula $C$, then the reflected copula $\nu_K(C)$ of $C$ is the distribution function of the random vector $ \boldsymbol{\eta}_{K} (\UUU,\eins - \UUU) $.
We refer to \cite{fuc2014gamma} for further details on the groups $\Gamma^{\pi}$, $\Gamma^{\nu}$ and $\Gamma$.



\subsection*{Concordance Order}

A relation $\lessdot$ on $\CC$ is said to be an \emph{order relation} if it is 
reflexive, antisymmetric and transitive;
in this case the pair $(\CC,\lessdot)$ is called \emph{ordered set}.
The copula $C \in \CC$ is said to be 
\begin{hylist}
\item
the \emph{greatest element}, if the inequality $D\lessdot C$ holds for all $D \in \CC$.

\item
the \emph{least element}, if the inequality $C\lessdot D$ holds for all $D \in \CC$.

\item
a \emph{maximal element}, if, for every $D \in \CC$, 
the inequality $C\lessdot D$ implies $C=D$.

\item
a \emph{minimal element}, if, for every $D \in \CC$, 
the inequality $D\lessdot C$ implies $C=D$.
\end{hylist}
We denote by
$$
  m(\CC,\lessdot)
$$
the set of all minimal elements of $(\CC,\lessdot)$.

\bigskip

For $ C,D \in \CC $ we write 
$ C \preceq D $ if $C(\uuu) \leq D(\uuu)$ and $(\tau(C))(\uuu) \leq (\tau(D))(\uuu)$ 
for every $ \uuu \in \I^d $.
Then $\preceq$ is an order relation which is called the \emph{concordance order} on $\CC$.
Since $\tau(M)=M$, the upper Fr{\'e}chet--Hoeffding bound $M$ is the greatest element 
in $(\mathcal{C}, \preceq)$; 
similarly, in the case $d=2$, the lower Fr{\'e}chet--Hoeffding bound $W$ is the least element 
in $(\mathcal{C}, \preceq)$. 
A copula $C \in \CC$ is said to be a 
\emph{minimal copula} if $C \in m(\CC,\preceq)$. 


\section{Multivariate Countermonotonicity}

\smallskip

With the absence of a least element for dimensions $d\geq3$, 
the set of all minimal elements, i.e. all locally least elements, in the collection of all copulas 
turns out to play an important role when studying extreme negative dependence concepts.
In this section, we discuss the existence and list some important examples of copulas that are minimal with respect to concordance order (so--called \emph{minimal copulas}). 
Additionally, we investigate several sufficient conditions and we provide a necessary condition for a copula to be minimal.

\bigskip

First sufficient and necessary conditions for a copula to be minimal can be achieved by comparing concordance order with other order relations:

\bigskip\pagebreak

\begin{remark}{}
Let $\lessdot$ be an order relation on $\CC$.
\begin{onelist}
\item
If $C \preceq D$ implies $C \lessdot D$, 
then 
$ m(\CC,\lessdot) \subseteq m(\CC,\preceq) $.

\item
If $C \lessdot D$ implies $C \preceq D$, 
then 
$ m(\CC,\preceq) \subseteq m(\CC,\lessdot) $.
\end{onelist}
Note that (1) is applicable to pointwise order and (2) is applicable to supermodular order; 
for more details on the comparison of concordance order with other order relations, we refer to \cite{joe1990, musc2000}.
\end{remark}{}

\bigskip

We are now interested in sufficient and necessary conditions for a copula to be minimal 
that are formulated in terms of the copula itself.

\bigskip

A copula $C \in\CC$ is said to be \emph{$K-$countermonotonic} ($K-$CM) 
if there exists 
\begin{hylist}
\item
some $ K \subseteq \{1,...,d\} $ with $ 2 \leq |K| \leq d $, 
\item
a family $\{g_k\}_{k \in K} $ of strictly increasing and continuous functions $ \I \to \R $ and 
\item
some $c \in \R$
\end{hylist}
such that 
$$
  Q^C \left[ \left\{ \uuu \in \I^d \, : \, \sum_{k \in K} g_k(u_k) = c \right\} \right]
	= 1
$$
compare \citep[Definitions 2,3 \& 4]{Ahn9}. 
We denote by $ \CC_{K-\rm CM} $ the collection of all copulas that are $K-$CM.
To ease notation, we write $d-$CM in the case $K=\{1,...,d\}$.
Note that, for $d=2$, a copula $C$ is $2-$CM if and only if $C=W$.
Thus, $K-$countermonotonicity may be regarded as a natural extension of countermonotonicity to dimensions $d\ge3$.

\bigskip

For some particular choices $K \subseteq \{1,...,d\}$, $K$--countermonotonicity implies minimality:

\bigskip

\begin{proposition}{} \label{dCMMinimal} \label{d-1CMMinimal} 
The inclusion
\begin{eqnarray*}
  \CC_{K-\rm CM}
	& \subseteq & m(\CC,\preceq) 
\end{eqnarray*}
holds for every $K \subseteq \{1,...,d\}$ such that $2 \leq |K| \in \{d-1,d\} $.
Moreover, 
$ \CC_{d-\rm CM} = m(\CC,\preceq) $
if and only if $d=2$.
\end{proposition}{}

\bigskip

\begin{proof}{}
The inclusions were proved in \citep[Theorem 4, Lemma 5]{Ahn9}, and
the equivalence for $d=2$ follows from Example \ref{dCMMinimalEx} (1) below.
\end{proof}{}

\bigskip

Examples \ref{dCMMinimalEx} below show that 
\begin{onelist}
\item
the inclusion in 
Proposition \ref{dCMMinimal} for $|K|=d$ is strict whenever $d \geq 3$, 

\item
the inclusion in 
Proposition \ref{dCMMinimal} for $|K|=d-1$ is strict, and that 

\item
the inclusion in Proposition \ref{dCMMinimal} fails to be satisfied whenever 
$ 2 \leq |K| \leq d-2 $.
\end{onelist}

\bigskip

In the following, we list some important minimal copulas and show in passing that the set $m(\CC, \preceq)$ is non--empty:

\bigskip

\begin{examples}{} \label{MinimalCopulasExample}
\begin{onelist} 
\item
For every $ K \subseteq \{1,...,d\} $ with $ 1 \leq |K| \leq d-1$, 
the copula 
$$ \nu_K(M) $$ is $d-$CM and hence minimal. 

Indeed, consider $ K \subseteq \{1,...,d\} $ with $ 1 \leq |K| \leq d-1$.
Then 
$$
  Q^{\nu_K(M)}\Bigl[ \Bigl\{ \uuu\in\I^d \colon \eeta_K(\uuu,\eins-\uuu) = \alpha\eins \;\text{for some   $ \alpha\in\I $} \Bigr\} \Bigr]  
  = 1           
$$
By choosing $ g_k(u_k) = u_k/|K| $ for every $ k\in K $ and $ g_k(u_k) = u_k/(d-|K|) $ for every $ k\in \{1,...,d\} \backslash K$,
we obtain
\begin{eqnarray*}
  \sum_{k=1}^{d} g_k(u_k) 
	& = & \sum_{k \in K} \frac{u_k}{|K|} + \sum_{k \in \{1,...,d\} \backslash K} \frac{u_k}{d-|K|}  
	\\
	& = & \sum_{k \in K} \frac{1-\alpha}{|K|} + \sum_{k \in \{1,...,d\} \backslash K} \frac{\alpha}{d-|K|} 
	\\
	& = & 1-\alpha + \alpha
	\\
	& = & 1
\end{eqnarray*}
for every $ \uuu \in \I^d $ satisfying $ \eeta_K(\uuu,\eins-\uuu) = \alpha\eins \;\text{for some   $ \alpha\in\I $}$, and hence
$$
  Q^{\nu_K(M)}\left[ \left\{ \uuu\in\I^d \colon \sum_{k=1}^{d} g_k(u_k) = 1 \right\} \right]
	= 1
$$
The assertion then follows from Proposition \ref{dCMMinimal}.

\item 
The copula   $ C $   given by 
$$
  C(\uuu)                                                             
  := \max \Biggl\{ \sum_{k=1}^d u_k^{1/(d-1)} - (d-1), 0 \Biggr\}^{d-1}  
$$
is Archimedean and 
is called the 
\emph{Clayton copula with parameter   $ -1/(d-1) $}.
It follows from \citep[Theorem 4]{Ahn7} and Proposition \ref{dCMMinimal} that $C$ is $d-$CM and hence minimal.

\item
Consider $d=3$ and the probability measure $Q: \BB(\I^3) \to \I$ whose probability mass is distributed uniformly on the edges of the equilateral triangle in $\I^3$ with vertices 
$(0,1/2,1)$, $(1/2,1,0)$ and $(1,0,1/2)$.
By \cite[Example 7]{neuf2012}, its corresponding distribution function is a copula satisfying
$$
  Q \left[ \left\{ \uuu \in \I^3 \, : \, \sum_{i=1}^{3} u_i = 3/2 \right\} \right]
	= 1
$$
which implies that $C$ is $3-$countermonotonic.
It hence follows from Proposition \ref{d-1CMMinimal} that $C$ is minimal.
\end{onelist}
\end{examples}{}

\bigskip

\begin{examples}{} \label{dCMMinimalEx}
\begin{onelist}
\item
Consider $d \geq 3$, 
a $(d-1)-$dimensional, $(d-1)-$CM copula $C$ and define the map $ D: \I^{d-1} \times \I \to \R $ by letting
$$
  D(\uuu,v)
	:= C(\uuu) \, v  
$$
Then, by \cite[Theorem 6.6.3]{Schweizer}, $D$ is a $d$--dimensional copula and 
it follows from Proposition \ref{d-1CMMinimal} that $D$ is minimal.
However, $D$ fails to be $d-$CM which follows from straightforward calculation.

\item
Consider $d \geq 3$.
Then the copula $C$ discussed in Example \ref{MinimalCopulasExample} (2) is $d-$CM and hence minimal.
However, it is evident that $C$ fails to be $K-$CM whenever $2 \leq |K| \leq d-1$.
\pagebreak

\item
Consider $d \geq 4$, $K \subseteq \{1,...,d\}$ with $2\leq|K|\leq d-2$,
a $|K|-$dimensional, $|K|-$CM copula $C$ and define the maps $ D,E: \I^{|K|} \times \I^{d-|K|} \to \R $ by letting
\begin{eqnarray*}
  D(\uuu,\vvv)
	& := & C(\uuu) \, \Pi(\vvv)  
	\\
	E(\uuu,\vvv)
	& := & C(\uuu) \, M(\vvv)
\end{eqnarray*}
Then, by \cite[Theorem 6.6.3]{Schweizer}, 
$D$ and $E$ are $d-$ dimensional copulas and, by definition, are $\{1,...,|K|\}-$CM.
However, $D \preceq E$ with $D \neq E$ which implies that $E$ is not minimal. 
\end{onelist}
\end{examples}{}

\bigskip

\begin{remark}{}
The collections $\CC_{K-\rm CM}$ with $K \subseteq \{1,...,d\}$ such that $2 \leq |K|$ 
are not directed.
Indeed, consider $K,L \subseteq \{1,...,d\} $ with $2 \leq |K|$ and $2 \leq |L|$. 
Then, due to Examples \ref{dCMMinimalEx} (1) and (2),
$K \subseteq L$ neither implies $\CC_{K-\rm CM} \subseteq \CC_{L-\rm CM}$ nor $\CC_{L-\rm CM} \subseteq \CC_{K-\rm CM}$.
\end{remark}{}

\bigskip

As shown in Example \ref{dCMMinimalEx} (3), 
the inclusion in Proposition \ref{dCMMinimal} fails to be satisfied whenever 
$ 2 \leq |K| \leq d-2 $.
Nevertheless, for every $d \geq 4$ and every $K \subseteq \{1,...,d\}$ with $ 2 \leq |K| \leq d-2 $, 
there exist $K-$CM copulas that are minimal.
The following construction principle generalizes Example \ref{dCMMinimalEx} (1):

\bigskip

\begin{theorem}{} \label{dCMProduct}
Consider $d \geq 4$, $K \subseteq \{1,...,d\}$ with $2 \leq |K| \leq d-2$, 
a $|K|-$dimensional, $|K|-$CM copula $C$ and a $(d-|K|)-$dimensional, $(d-|K|)-$CM copula $D$.
Then the copula $E: \I^{|K|} \times \I^{d-|K|}$ given by 
$$
  E(\uuu,\vvv) 
	:= C(\uuu) \, D(\vvv)
$$
is $d-$CM and hence minimal.
\end{theorem}{}

\bigskip

\begin{proof}{}
By, \cite[Theorem 6.6.3]{Schweizer}, $E$ is a copula and,
by definition,
there exist families 
$\{g_k\}_{k \in \{1,...,|K|\}}$ and 
$\{h_l\}_{l \in \{1,...,d-|K|\}}$ of strictly increasing and continuous functions $ \I \to \R $ \linebreak
and 
constants $c,d \in \R$ such that 
$ Q^C \big[ \big\{ \uuu \in \I^{|K|} \, : \, \sum_{k \in \{1,...,|K|\}} g_k(u_k) = c \big\} \big]
	= 1 $
and \linebreak
$ Q^D \big[ \big\{ \uuu \in \I^{d-|K|} \, : \, \sum_{l \in \{1,...,d-|K|\}} h_l(u_l) = d \big\} \big]
	= 1 $.
Then the copula measure $Q^E$ of $E$ satisfies
\begin{eqnarray*}
  1
	& \geq & Q^E \biggl[ \biggl\{ \uuu \in \I^d \, : \, \sum_{k \in \{1,...,|K|\}} g_k(u_k) 
					 + \sum_{l \in \{1,...,d-|K|\}} h_l(u_{|K|+l}) = c+d \biggr\} \biggr]
  \\
	& \geq & Q^C \biggl[ \biggl\{ \uuu \in \I^{|K|} \, : \, \sum_{k \in \{1,...,|K|\}} g_k(u_k) = c\biggr\} \biggr] 
	         \; Q^D \biggl[ \biggl\{ \uuu \in \I^{d-|K|} \, : \, \sum_{l \in \{1,...,d-|K|\}} h_l(u_l) = d \biggr\} \biggr]
	\\
	&   =  & 1
\end{eqnarray*}
Thus, $E$ is $d-$CM and it follows from Proposition \ref{dCMMinimal} that $E$ is minimal.
\end{proof}{}

\bigskip

We proceed with the discussion of a necessary condition for a copula to be minimal.

\bigskip

A copula $C\in\CC$ is said to be \emph{Kendall-countermonotonic} ($\tau$-CM) if the identity
$$
  \min \big\{ C(\uuu), (\tau(C)) ({\bf 1}-\uuu) \big\} = 0
$$
holds for every $ \uuu \in ({\bf 0},{\bf 1}) $.
We denote by $ \CC_{\tau-{\rm CM}} $ the collection of all copulas that are $\tau$-CM.
The term Kendall-countermonotonicity is motivated by the fact that a copula $C$ is $\tau$-CM if and only if $C$ minimizes multivariate Kendall's tau; see \citep{fumcsc2018}.  \pagebreak 
Kendall's tau is a map $\kappa: \CC \to \R$ given by
\begin{eqnarray*}
  \kappa(C)                                                       
  & := &  \frac{2^{d}}{2^{d-1}-1}\,\biggl( \int_{\I^d} C(\uuu) \; \mathrm{d} Q^{C} (\uuu) - \frac{1}{2^d} \biggr)  
\end{eqnarray*}
and the definition of Kendall's tau is in accordance with \cite{nel2002}. 
The following characterization of Kendall-countermonotonicity is due to \citep[Lemma 3.1, Theorem 3.4]{fumcsc2018}:

\bigskip

\begin{proposition}{} \label{MinimalCopulas}
For a copula $C \in \CC$ the following are equivalent:
\begin{onelist}
\item
$C$ is $\tau$--CM.

\item
$\tau(C)$ is $\tau$--CM.

\item
Every $ \uuu \in ({\bf 0},{\bf 1}) $ satisfies
$ Q^C [[{\bf 0}, \uuu]] = 0 $ or $	Q^C [[\uuu,{\bf 1}]] = 0 $.

\item 
$\int_{\I^d} C(\uuu) \; \mathrm{d} Q^{C} (\uuu) = 0$.
\end{onelist}
\end{proposition}{}

\bigskip

It follows from \citep[Theorem 3.3]{fumcsc2018} that, for $d=2$, 
a copula $C$ is $\tau$-CM if and only if $C=W$.
Thus, Kendall-countermonotonicity may be regarded as a natural extension of countermonotonicity to dimensions $d\ge3$.

\bigskip

\begin{remark}{} \label{KendallRemark}
A subset $A \subseteq \I^d$ is said to be \emph{strictly comonotonic} 
if either $\uuu < \vvv$ or $\vvv < \uuu$ for all $\uuu,\vvv \in A$.
Consider now a copula $C$ for which there exists some strictly comonotonic set 
$A \subseteq ({\bf 0},{\bf 1})$ consisting of at least two points such that 
$$
  A \subseteq {\rm supp} \; Q^C 
$$ 
i.e. the support ${\rm supp} \; Q^C $ of $Q^C$ contains some strictly comonotonic subset of $({\bf 0},{\bf 1})$.
It then follows from Proposition \ref{MinimalCopulas} (3) that such a copula
$C$ fails to be Kendall--countermonotonic. 
\end{remark}{}

\bigskip

Thus, one may interpret Kendall--countermonotonicity as the one extreme negative dependence concept 
where it is inadmissible for a copula to have some strictly comonotonic support.

\bigskip

The following theorem is the main result of this paper; 
it states that every minimal copula is $\tau$-CM and hence every minimal copula minimizes Kendall's tau.

\bigskip

\begin{theorem}{} \label{MinimalKendallCM}
We have 
\begin{eqnarray*}
  m(\CC,\preceq) 
	& \subseteq & \CC_{\tau-{\rm CM}} 
\end{eqnarray*}
Moreover, if $d=2$, then $ m(\CC,\preceq) = \CC_{\tau-{\rm CM}} $.
\end{theorem}

\bigskip

\begin{proof}
The inclusion is proved in the appendix (see Lemmas \ref{A1}, \ref{A2} and \ref{A3});
there, for a copula $ C \in \CC \backslash \CC_{\tau-\rm CM} $, 
we construct a copula $D \in \CC$ satisfying $ D \preceq C $ with $D \neq C$ which then implies that $C$ is not minimal. 
The identity for $d=2$ follows from \citep[Theorem 3.3]{fumcsc2018}.
\end{proof}

\bigskip

The following example shows that the inclusion in 
Theorem \ref{MinimalKendallCM} is strict whenever $d \geq 4$:

\bigskip

\begin{example}{} \label{MinimalKendallCMEx}
For $d \geq 4$, consider the maps $ C,D: \I^2 \times \I^{d-2} \to \I $ given by 
\begin{eqnarray*}
  C(\uuu,\vvv)
	& := & W(\uuu) \, \Pi(\vvv)
	\\* 
	D(\uuu,\vvv)
	& := & W(\uuu) \, M(\vvv)
\end{eqnarray*}
By, \cite[Theorem 6.6.3]{Schweizer}, $C$ and $D$ are copulas,
and it follows from \citep[Remark 3.1.(3)]{fumcsc2018} that $D$ is $\tau$--CM.
However, $D$ fails to be a minimal copula which is a consequence of $ C \preceq D $ with $C \neq D$. 
\end{example}{} \pagebreak

\bigskip

It is interesting to note that also every $K-$CM copula is $\tau-$CM:

\bigskip

\begin{proposition}{} \label{KCMTauCM}
The inclusion
\begin{eqnarray*}
  \CC_{K-{\rm CM}} 
	& \subseteq & \CC_{\tau-{\rm CM}} 
\end{eqnarray*}
holds for every $ K \subseteq \{1,...,d\} $ such that $ 2 \leq |K| \leq d $.
Moreover, $ \bigcup_{K \subseteq \{1,...,d\}, 2 \leq |K| \leq d} \; \CC_{K-{\rm CM}} = \CC_{\tau-{\rm CM}} $ if and only if $d=2$.
\end{proposition}{}

\bigskip

\begin{proof}{}
The inclusions follow from \cite[Theorem 6]{Ahn7},
and the equivalence for $d=2$ is a consequence of Example \ref{KCMTauCMEx} below.
\end{proof}

\bigskip

The following example shows that the inclusions in 
Proposition \ref{KCMTauCM} are strict whenever $d \geq 3$:

\bigskip

\begin{example}{} \label{KCMTauCMEx}
For $d \geq 3$, consider the copula
$$
  C(\uuu)
	:= \frac{1}{2^d-2} \; \sum_{K \subseteq \{1,...,d\}, 1 \leq |K| \leq d-1} \nu_K (M) 
$$
Then, $C$ is $\tau$--CM 
but fails to be $K$--CM whenever $ 2 \leq |K| \leq d $.
The result is a $d-$dimensional analogue of \cite[Example 3.2]{fumcsc2018}.
\end{example}{}

\bigskip

In summary, we thus have

\bigskip

\begin{center}
\begin{tabular}{ccc}
  $m(\CC,\preceq)$ 
	& $\subseteq_{II}$
	& $\CC_{\tau-{\rm CM}}$
	\\ \\
	$\rotatebox{90}{$\subseteq$}_{I}$ 
	&& $\rotatebox{90}{$\subseteq$}_{IV}$
	\\ \\
	$\bigcup_{K \subseteq \{1,...,d\}, 2 \leq |K| \in \{d,d-1\}} \; \CC_{K-{\rm CM}}$
	& $\quad \subseteq_{III} \quad$ 
	& $\bigcup_{K \subseteq \{1,...,d\}, 2 \leq |K| \leq d} \; \CC_{K-{\rm CM}}$
\end{tabular}
\end{center}

\bigskip

Recall that, for $d=2$, all the sets are identical and that,
for $d \geq 4$, inclusions $II$, $III$ and $IV$ are strict. 
In the case $d=3$, the above figure reduces to   
$$
  \bigcup_{K \subseteq \{1,...,d\}, 2 \leq |K| \leq 3} \; \CC_{K-{\rm CM}}
	\qquad \subseteq \qquad m(\CC,\preceq)
	\qquad \subseteq \qquad \CC_{\tau-{\rm CM}}
$$
such that, due to Example \ref{KCMTauCMEx}, at least one of the inclusions is strict.

\section{Continuous and Order Preserving Functionals}

\smallskip

In this section we study minimal copulas in connection with continuous and concordance order preserving functionals and show that each such optimization problem is minimized by some minimal copula.
In particular, we show that any continuous functional that is even strictly concordance order preserving is minimized by minimal copulas only.

\bigskip \pagebreak

A map $ \kappa: \CC \to \R $ is said to be 
\begin{hylist} 
\item 
\emph{continuous} if, for any sequence $ \{C_n\}_{n \in \N} \subseteq \CC $ and any copula $C \in \CC$,
uniform convergence $ \lim_{n \to \infty} C_n = C $ implies 
$ \lim_{n \to \infty} \kappa(C_n) = \kappa(C) $.


\item
\emph{concordance order preserving} if the inequality $ \kappa(C) \leq \kappa(D) $ holds for all 
$ C,D \in \mathcal{C} $ satisfying $ C \preceq D $.


\item 
\emph{strictly concordance order preserving} if it is concordance order preserving and the strict inequality $ \kappa(C) < \kappa(D) $ holds for all 
$ C,D \in \mathcal{C} $ satisfying $ C \preceq D $ with $C \neq D$.
\end{hylist}

\bigskip

We start with the discussion of some topological properties of $\CC$:
It is well--known that
$\CC$ is a compact subset of the space $ (\Xi(\I^d), d_{\infty})$ of all continuous real--valued functions with domain $\I^d$ under the topology of 
uniform convergence; 
see, e.g., 
\citep[Theorem 1.7.7]{dus2016}. 
It is further well--known that 
the range of $\CC$ with respect to any continuous map $ \kappa: \CC \mapsto \R $ is compact in $\R$;
see, e.g., \citep[Theorem 4.14]{Rudin1976}.

For a given map $ \kappa: \CC \mapsto \R$ and a subset $ \DD \subseteq \CC $, we define the set 
$ m (\kappa, \DD) $ by letting
$$
  m(\kappa, \DD)
	:= \left\{ D \in \DD \, : \, \kappa(D) = \inf \limits_{C\in \DD} \kappa(C) \right\}
$$
Note that every continuous functional $\CC \mapsto \Real$ is minimized by some copula $C \in \CC$;
see, e.g., \citep[Theorems 4.15 and 4.16]{Rudin1976}:

\bigskip

\begin{proposition}{} \label{comp.thm}
Let $\kappa: \CC \mapsto \Real$ be continuous and $\DD \subseteq \CC$ be a compact subset of $ (\CC, d_{\infty}) $.\linebreak
Then 
$ m(\kappa, \DD)$ is non--empty. 
\end{proposition}{}



\bigskip

Even though a continuous functional $\CC \to \R$ is minimized 
by some copula $C \in \CC$,
the calculation of its minimal value can be quite difficult. 
For an illustration, let us consider the following quite popular optimization problem for which a solution was recently presented in \cite{wawa2011}:

\bigskip

\begin{example}{} \label{ExPreSpearman}
The map $ \kappa: \CC \to \R $ given by
$$
 \kappa (C) 
  := \int_{\I^d} \Pi(\uuu) \; \mathrm{d} Q^C (\uuu) 
$$
is a continuous and strictly concordance order preserving functional,
and thus, due to Proposition \ref{comp.thm}, $\kappa$ is minimized by some copula $C \in \CC$.
In \cite[Corollary 4.1 and Figure 3.2]{wawa2011} the authors have presented a solution for
$ \inf_{C \in \CC} \kappa(C) $ in arbitrary dimension
and, for $d=3$, a minimal copula minimizing $\kappa$.
\end{example}



\bigskip

We now show that, for any compact subset $\DD \subseteq \CC $, 
the ordered set $(\DD, \preceq)$ is \emph{coverable from below}, i.e. 
for every copula $D \in \DD$, there exists some minimal copula $C \in (\DD, \preceq)$ satisfying
$ C \preceq D $.

\bigskip

\begin{theorem}{}\label{existence.thm}
If $ \DD \subseteq \CC $ is compact, 
then $ (\DD, \preceq) $ is coverable from below.
In particular, $ (\CC, \preceq) $ is coverable from below.
\end{theorem}

\bigskip

\begin{proof}
Consider $ D \in \DD $ and define the subset $ \EE \subseteq \DD $ by letting
$ \EE
	:= \left\{E \in \DD \big\vert E \preceq D \right\} $.
Since $ D \in \EE $, the set $\EE$ is non--empty.
We first show that $\EE$ is a compact subset of $ (\DD, d_{\infty}) $.

Since $\EE$ is a subset of the compact set $\DD$, \pagebreak
it is enough to show that
$\EE$ is closed; compare \citep[Theorem 2.35]{Rudin1976}.
To this end, consider a convergent sequence of copulas $ \{E_n\}_{n \in \N} \subseteq \EE $ 
and define the map $ E_\infty: \I^d \to \R $ by letting
$ E_\infty (\uuu)
	:= \lim_{n\rightarrow \infty} E_n (\uuu) $.
Then, by \citep[Theorem 1.7.5]{dus2016}, 
$ E_\infty $ is a copula,
by \citep[Theorem 1.7.6]{dus2016}, $E_\infty \in \DD$,
and, by definition, $E_\infty$ satisfies
$ E_\infty \preceq D $.
Thus, $E_\infty \in \EE$, and 
\citep[Theorem 1.7.6]{dus2016} implies that $\EE$ is closed and hence compact.

Now, consider the continuous and strictly concordance order preserving functional $\kappa$ discussed in Example \ref{ExPreSpearman}.
Then, by compactness of $\EE$ and Proposition \ref{comp.thm}, 
there exists some $E_0\in \EE \subseteq \DD$ satisfying
$ \kappa (E_0)
	= \inf_{E \in \EE} \kappa (E) $.
To show that $E_0$ is a minimal copula in $ (\DD, \preceq)$, 
consider some copula $ E_1 \in \DD $ satisfying $ E_1 \preceq E_0 $.
Then $ E_1 \in \EE $.
Assuming $ E_1 \neq E_0 $, 
the fact that $\kappa$ is strictly concordance order preserving implies 
$ \kappa (E_1) < \kappa (E_0) $
which contradicts $\kappa (E_0) = \inf_{E \in \EE} \kappa (E)$. 
So we conclude that $E_1=E_0$ which implies that $E_0$ is a minimal copula in $(\DD, \preceq)$ and satisfies $ E_0 \preceq D $.
This proves the assertion.
\end{proof}

\bigskip

\begin{remark}{}
For a minimal copula $D \in m(\CC,\preceq)$, define the set 
$\CC_D:=\{ C \in \CC : D \preceq C \}$.
Then, Theorem \ref{existence.thm} yields
$$
  \CC = \bigcup_{D \in m(\CC,\preceq)} \CC_D 
$$
i.e. every copula $C \in \CC$ is comparable to (at least) one minimal copula $D \in m(\CC,\preceq)$.
\end{remark}{}

\bigskip

The following result shows that the set of minimal copulas plays a key role when searching for the minimal value of a continuous and (strictly) concordance order preserving functional.
It turns out that, for any continuous and concordance order preserving map $\kappa$,
the set of minimizers contains at least one minimal copula.
It further turns out that any continuous map that is even strictly concordance order preserving is minimized by minimal copulas only.

\bigskip

\begin{theorem}{} \label{FunctionalCOP.1} \label{FunctionalCOP.2}
Let $ \kappa: \CC \mapsto \R $ be a continuous and concordance order preserving map.
\begin{1onelist}
\item 
If $ \DD \subseteq \CC $ is compact, 
then $m(\kappa, \DD)$ contains at least one minimal copula of $ (\DD, \preceq) $.

\item
In particular, 
$m(\kappa, \CC)$ contains at least one minimal copula of $ (\CC, \preceq) $.

\item 
If $d=2$, then
$ \{W\} = m(\CC,\preceq) \subseteq m(\kappa, \CC) $. 
\end{1onelist}
Further assume that $\kappa$ is also strictly concordance order preserving.
\begin{2onelist}
\item 
If $ \DD \subseteq \CC $ is compact, 
then $ m(\kappa,\DD) $ is a non--empty set of minimal copulas of $(\DD, \preceq)$.

\item 
In particular, 
$ m(\kappa,\CC) $ is a non--empty set of minimal copulas of $(\CC, \preceq)$.
	
\item 
If $d=2$, then 
$ \{W\} = m(\CC,\preceq) = m(\kappa, \CC) $. 
\end{2onelist}
\end{theorem}{}

\bigskip

\begin{proof}
The first three assertions follow from Proposition \ref{comp.thm} and Theorem \ref{existence.thm}.

We now prove the second part. To this end, assume that $\kappa$ is strictly concordance order preserving.
Since $m(\kappa,\DD)$ is non--empty, 
by Proposition \ref{comp.thm}, there exists some copula $ D_0 \in m(\kappa,\DD) $ satisfying 
$ \kappa(D_0)
	= \inf_{D \in \DD} \kappa(D) $.
To show that $D_0$ is a minimal copula in $ (\DD, \preceq)$, consider some copula $D_1 \in \DD$ satisfying $ D_1 \preceq D_0 $.
Then $ D_1 \in m(\kappa,\DD) $.
Further, assume that $ D_1 \neq D_0 $.
Since $\kappa$ is strictly concordance order preserving we hence obtain
$ \kappa(D_1)
	< \kappa(D_0) $
which contradicts $\kappa(D_0) = \inf_{D \in \DD} \kappa(D)$.
Therefore, $D_1=D_0$, which concludes that $D_0$ is a minimal copula in $(\DD, \preceq)$.
This proves the assertion.
\end{proof}

\bigskip

Note that the results (2.1) and (2.2) in Theorem \ref{FunctionalCOP.2} can not be extended to concordance order preserving 
functionals that fail to be strictly concordance order preserving 
since a copula attaining the minimal value of such a functional may fail to be a minimal copula; 
see, e.g., Corollary \ref{MOCKendallStrictMOC} together with Example \ref{MinimalKendallCMEx}.

Further note that, for $d \geq 3$, the inclusion in Theorem \ref{FunctionalCOP.2} (2.2) is strict, in general; 
see, e.g., Example \ref{MOCMinimalKendall.Ex}.


\section{Measures of Concordance} \label{SectionMOC}

\smallskip

In the following we apply the results of the previous sections to measures of concordance:
First of all, we show that in the class of all continuous and concordance order preserving measures of concordance Kendall's tau is particular since it is minimized by every minimal copula.
As a consequence of Theorem \ref{FunctionalCOP.1}, it further turns out that every continuous and strictly concordance order preserving measure of concordance is minimized by minimal copulas only.
Since the latter result is applicable to Spearman's rho we may conclude that every copula minimizing Spearman's rho is also a minimizer of Kendall's tau.

\bigskip

We employ the quite general definition of a measure of concordance proposed in \citep{fuc2016moc}; compare also \citep{duf2006, tay2007, tay2016}:
A map
$ \kappa: \CC \to \R $
is said to be a {\it measure of concordance} \index{measure of concordance} if it satisfies the following axioms:
\begin{romlist}
\item 
$ \kappa(M) = 1 $.

\item 
The identity 
$ \kappa (\gamma(C)) = \kappa (C) $
holds for all $ \gamma \in \Gamma^{\pi} $ and all $ C \in \CC $.

\item 
The identity 
$ \kappa (\tau(C)) = \kappa (C) $
holds for all $ C \in \CC $.

\item 
The identity 
$ \sum_{\nu \in \Gamma^{\nu}} \kappa (\nu(C)) = 0 $
holds for all $ C \in \CC $.  
\end{romlist} 
For the case $d=2$, this definition is in accordance with 
\citep{emt2004, emt2005, fus2014, sca1984}.



\bigskip

\begin{example}{} {\rm (Kendall's tau)} \label{ExampleKendallDef} \quad
The map $ \kappa^{(\tau)}: \CC \to \R $ given by 
$$
  \kappa^{(\tau)}(C) 
	:= \frac{2^d}{2^{d-1}-1} \left( \int_{\I^d} C(\uuu) \; \mathrm{d} Q^C(\uuu) - \frac{1}{2^d} \right)
$$
is a continuous and concordance order preserving measure of concordance, 
and is called \emph{Kendall's tau}; 
this definition of Kendall's tau is in accordance with that in \cite{nel2002}. 
$\kappa^{(\tau)}$ satisfies
$$
  \min_{C \in \CC} \kappa^{(\tau)}(C) = - \, \frac{1}{2^{d-1}-1}
$$
and its minimum value is attained by the minimal copula $\nu_1(M)$; 
see \cite[Theorem 5.1]{ubf2005} and \cite[Remark 3.1 (1)]{fumcsc2018}.
It follows from Example \ref{MOC.Ex1} that Kendall's tau is not strictly concordance order preserving.
\end{example}

\bigskip

\begin{example}{} {\rm (Spearman's rho)} \label{ExampleSpearmanDef} \quad
The map $ \kappa^{(\rho)}: \CC \to \R $ given by
$$
 \kappa^{(\rho)} (C) 
  := \frac{2^d\,(d+1)}{2^d-(d+1)} \left( \int_{\I^d} \frac{C(\uuu)+(\tau(C))(\uuu)}{2}  \; \mathrm{d} Q^{\Pi}(\uuu) - \frac{1}{2^d} \right)
$$
is a continuous and strictly concordance order preserving measure of concordance; 
the definition of Spearman's rho used here is in accordance with that in \cite{nel2002}.
Even though $\kappa^{(\rho)}$ can be minimized only by minimal copulas 
(Theorem \ref{FunctionalCOP.2}), its minimal value is known only for $d \in \{2,3\}$: 
In the case $d=2$, $W$ is the only copula minimizing Spearman's rho with $\kappa^{(\rho)} (W)=-1$ and, 
for $d=3$, \citet[Theorem 4]{neuf2012} have shown that 
$$
  \min_{C \in \CC} \kappa^{(\rho)} (C) = - \frac{1}{2}
$$
and that this minimal value is attained only by those minimal copulas satisfying \linebreak
$ Q [\{ \uuu \in \I^3 \, : \, \sum_{i=1}^{3} u_i = 3/2 \}] = 1$.
\end{example} \pagebreak

\bigskip

It immediately follows from Theorem \ref{FunctionalCOP.2} and Theorem \ref{MinimalKendallCM} that every continuous and strictly concordance order preserving measure of concordance is minimized by minimal copulas only,
and that every minimal copula minimizes Kendall's tau:

\bigskip

\begin{corollary}{} \label{MOCKendallStrictMOC}
\begin{onelist}
\item
The inclusions
$$
  m(\kappa,\CC) 
	\quad \subseteq \quad m(\CC,\preceq)
	\quad \subseteq \quad m(\kappa^{(\tau)},\CC)
$$
hold for every continuous and strictly concordance order preserving measure of concordance $\kappa$.

\item 
If $d=2$, then the identities
$$ 
  m(\kappa,\CC) 
	\quad = \quad m(\CC,\preceq) 
	\quad = \quad m(\kappa^{(\tau)},\CC) 
$$ 
hold for every continuous and strictly concordance order preserving measure of concordance $\kappa$.
\end{onelist}
\end{corollary}{}

\bigskip

It follows from Example \ref{MinimalKendallCMEx} 
that the second inclusion in Corollary \ref{MOCKendallStrictMOC} (1) is strict whenever $d \geq 4$. 
The next example shows that the first inclusion in Corollary \ref{MOCKendallStrictMOC} (1) is strict whenever $d \geq 3$:

\bigskip

\begin{example}{} \label{MOCMinimalKendall.Ex}
For every $d \geq 3$, 
there exist minimal copulas $ C,D \in m(\CC,\preceq)$ satisfying
$ \kappa^{(\rho)} (C) < \kappa^{(\rho)} (D) $. In particular,
$$
   m(\kappa^{(\rho)},\CC) \neq m(\CC,\preceq)
$$
Indeed, for $d=3$, consider the minimal copula $\nu_1(M)$ and the minimal copula $C$ discussed in  
Example \ref{MinimalCopulasExample} (3).
Then, by \cite[Example 7]{neuf2012}, we obtain
$$
  \kappa^{(\rho)} (C)
	= - \, \frac{1}{2}
	< - \, \frac{1}{3} 
	= \kappa^{(\rho)} (\nu_1(M)) 
$$
For $d \geq 4$, the minimal copulas $\nu_1(M)$ and $\nu_{1,2}(M)$ satisfy
$$
  \kappa^{(\rho)} (\nu_{1,2}(M))
	= \frac{2^{d+1}-(d-1)d(d+1)}{(d-1)d \, (2^d-(d+1))}
  < \frac{2^d-d(d+1)}{d \, (2^d-(d+1))} 
  = \kappa^{(\rho)} (\nu_1(M)) 
$$
This proves the assertion.
\end{example}{}

\bigskip

The relationship between measures of concordance, 
in particular between bivariate Kendall's tau and bivariate Spearman's rho,
has received considerable attention in literature; see, e.g., \cite{Caperaa, Fredricks, Hurlimann, schpautru2017}.
We are able to contribute to this topic by showing that every minimizer of Spearman's rho is also a minimizer of Kendall's tau:

\bigskip

\begin{corollary}{} \label{MOCKendallSpearman}
We have
$$ m(\kappa^{(\rho)},\CC) 
	\quad \subseteq \quad m(\kappa^{(\tau)},\CC) 
$$
Moreover, 
$ m(\kappa^{(\rho)},\CC) 
	= m(\kappa^{(\tau)},\CC) $ if and only if $d=2$.
\end{corollary}{}

\bigskip

For $d \in \{2,3\}$, the results in Corollary \ref{MOCKendallSpearman} are in accordance with 
Examples \ref{ExampleKendallDef} and \ref{ExampleSpearmanDef}.

\bigskip\pagebreak

We conclude this section by showing that Kendall's tau is not strictly concordance preserving:

\bigskip

\begin{example}{} \label{MOC.Ex1}\label{MOC.Ex2}
Consider the copula $ A: \I^2\to\I $ 
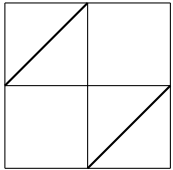
\begin{center}
\begin{tikzpicture}[xscale=0.55,yscale=0.55] 
\draw [thin]  (0,0) -- (0,4);
\draw [thin]  (2,0) -- (2,4);
\draw [thin]  (4,0) -- (4,4);
\draw [thin]  (0,0) -- (4,0);
\draw [thin]  (0,2) -- (4,2);
\draw [thin]  (0,4) -- (4,4);
\draw [thick] (0,2) -- (2,4);
\draw [thick] (2,0) -- (4,2);
\end{tikzpicture}
\end{center}
defined as the shuffle of   $ M $  with respect to the shuffling structure   
$ \{[\aaa_i,\bbb_i]\}_{i\in\{1,2\}} $   with 
$$\begin{array}{ccccccc}
 \aaa_1 &=& (\phantom{00}0,          1/2)  &&  \bbb_1 &=& (1/2,\phantom{00}1)  \\
 \aaa_2 &=& (1/2          ,\phantom{00}0)  &&  \bbb_2 &=& (\phantom{00}1,1/2)  
\end{array}$$
and the copula   $ B : \I^2\to\I $   
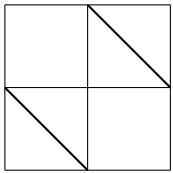
\begin{center}
\begin{tikzpicture}[xscale=0.55,yscale=0.55] 
\draw [thin]  (0,0) -- (0,4);
\draw [thin]  (2,0) -- (2,4);
\draw [thin]  (4,0) -- (4,4);
\draw [thin]  (0,0) -- (4,0);
\draw [thin]  (0,2) -- (4,2);
\draw [thin]  (0,4) -- (4,4);
\draw [thick] (0,2) -- (2,0);
\draw [thick] (2,4) -- (4,2);
\end{tikzpicture}
\end{center}
defined as the shuffle of $ W $ with respect to the shuffling structure   
$ \{[\aaa_i,\bbb_i]\}_{i\in\{1,2\}} $   with 
$$\begin{array}{ccccccc}
 \aaa_1 &=& (\phantom{00}0,\phantom{00}0)  &&  \bbb_1 &=& (1/2,1/2)  \\
 \aaa_2 &=& (1/2,1/2)  &&  \bbb_2 &=& (\phantom{00}1,\phantom{00}1)  
\end{array}$$
For more details on shuffles of copulas we refer to \cite{dufs2010}.
The copulas $A$ and $B$ satisfy \linebreak $A(\uuu) \leq B(\uuu)$ for all $\uuu \in \I^2$ and hence 
$A \preceq B$ with $A \neq B$ and, 
by \cite[Corollary 5.2]{fumcsc2018}, we obtain 
\begin{eqnarray*}
  \kappa^{(\tau)}(A) 
	& = & \kappa^{(\tau)}(B) 
\end{eqnarray*}
Moreover, for $d \geq 3$, define the functions $C, D: \I^2 \times \I^{d-2} \to \I $ by letting
\begin{eqnarray*}
  C(\uuu,\vvv) 
	& := & A(\uuu) \, \prod_{i=1}^{d-2} v_i
	\\
	D(\uuu,\vvv) 
	& := & B(\uuu) \, \prod_{i=1}^{d-2} v_i
\end{eqnarray*}
By \cite[Theorem 6.6.3]{Schweizer}, $C$ and $D$ are copulas,
$C \preceq D$ with $C \neq D$, and, by \cite[Corollary 5.2]{fumcsc2018}, we obtain
\begin{eqnarray*}
  \kappa^{(\tau)}(C) 
	& = & \kappa^{(\tau)}(D) 
\end{eqnarray*}
Thus, Kendall's tau is not strictly concordance order preserving.

\end{example}{}

\bigskip

Note that, for $d \geq 4$, the non--strictness of Kendall's tau also follows from 
Theorem \ref{FunctionalCOP.1} and
Example \ref{MinimalKendallCMEx}.

\appendix
\section{Appendix}

\smallskip

In this section we prove Theorem \ref{MinimalKendallCM}; 
the idea of the proof is as follows:
For a copula $ C \in \CC \backslash \CC_{\tau-\rm CM} $,  
we construct a copula $D \in \CC$ satisfying $ D \preceq C $ with $D \neq C$ which then implies that $C$ is not a minimal copula. 
More precisely, we extract some comonotonic part of the given copula $C$ 
and construct a related copula $D$ in which this comonotonic part is "made non--comonotonic".

\bigskip

The following result provides the basis for the construction of the copula $D$:

\bigskip

\begin{lemma}{} \label{A1}
For every copula $ C \in \CC \backslash \CC_{\tau-\rm CM} $ 
there exist some $ p \in (0,0.5] $ and 
some $\aaa, \bbb \in ({\bf 0}, {\bf 1}) $ with $ \aaa \leq \bbb $ such that
$ Q^C [[{\bf 0},\aaa]]
	= p 
	= Q^C [[\bbb, {\bf 1}]] $.
\end{lemma}{}

\bigskip

\begin{proof}{}
Consider $ C \in \CC \backslash \CC_{\tau-\rm CM} $.
Then, by Proposition \ref{MinimalCopulas},
there exists some $\uuu \in ({\bf 0}, {\bf 1})$ satisfying 
$ C(\uuu) > 0 $ and $ Q^C [[\uuu,{\bf 1}]] > 0 $.
W.l.o.G, let $ Q^C [[\uuu,{\bf 1}]] \leq C(\uuu) $ and put $ p:= Q^C [[\uuu,{\bf 1}]] $.
Then $ p \in (0,0.5] $.
Since every copula is continuous, the map $ \I \to \I $ given by
$ \alpha \mapsto C (\alpha \uuu) $
is continuous as well.
Thus, there exists some $ \beta \in (0,1] $ satisfying 
$ \beta \uuu \in ({\bf 0}, {\bf 1}) $,
$ \beta \uuu \leq \uuu $ and
$ C (\beta \uuu) = p $.
This proves the assertion.
\end{proof}{}

\bigskip

For a copula $ C \in \CC \backslash \CC_{\tau-\rm CM} $ satisfying
$ Q^C [[{\bf 0},\aaa]]
	= p 
	= Q^C [[\bbb, {\bf 1}]] $
for some $ p \in (0,0.5] $ and some $\aaa, \bbb \in ({\bf 0}, {\bf 1}) $ with $ \aaa \leq \bbb $,
we first define the maps $C_{(\aaa)}$, $C_{(\bbb)}: \I^d \to \I $ by letting
\begin{eqnarray*}
  C_{(\aaa)} (\uuu)
	& := & \frac{1}{p} \; Q^C \big[ [{\bf 0}, \uuu] \cap [{\bf 0}, \aaa] \big]
	\\
	C_{(\bbb)} (\uuu)
	& := & \frac{1}{p} \; Q^C \big[ [{\bf 0}, \uuu] \cap [\bbb, {\bf 1}] \big] 
\end{eqnarray*}
Then $C_{(\aaa)}$ and $C_{(\bbb)}$ are $d$--dimensional distribution functions on $\I^d$.
We further define the maps $C_{(1,\aaa,\bbb)}$, $C_{(2,\aaa,\bbb)}: \I^d \to \I $ by letting
\begin{eqnarray*}
  C_{(1,\aaa,\bbb)} 
	& := & \frac{1}{2} \; C_{(\aaa)} + \frac{1}{2} \; C_{(\bbb)}
\end{eqnarray*}
and
\begin{eqnarray*}
	C_{(2,\aaa,\bbb)} (\uuu)
	& := & \frac{1}{2} \;  
				       C_{(\aaa)} \big( \boldsymbol{\eta}_1 ({\bf 1},\uuu) \big) 
	       \cdot C_{(\bbb)} \big( \boldsymbol{\eta}_1 (\uuu,{\bf 1}) \big) 
				   +   \frac{1}{2} \;
					     C_{(\bbb)} \big( \boldsymbol{\eta}_1 ({\bf 1},\uuu) \big) 
				 \cdot C_{(\aaa)} \big( \boldsymbol{\eta}_1 (\uuu,{\bf 1}) \big) 
	\\*
	&  = & \frac{1}{2} \; \bigg( \frac{1}{p} \;  Q^C \big[ [0, u_1] \times \I^{d-1} \cap [{\bf 0}, \aaa] \big] \cdot
				 \frac{1}{p} \;  Q^C \big[ \I \times [0,u_2] \times [0,u_d] \cap [\bbb,{\bf 1}] \big]
	\\*
	&    & \quad + \frac{1}{p} \;  Q^C \big[ [0, u_1] \times \I^{d-1} \cap [\bbb,{\bf 1}] \big] \cdot
				 \frac{1}{p} \;  Q^C \big[ \I \times [0,u_2] \times [0,u_d] \cap [{\bf 0}, \aaa] \big] \bigg)
\end{eqnarray*}
Then $C_{(1,\aaa,\bbb)}$ and $C_{(2,\aaa,\bbb)}$ are also $d$--dimensional distribution functions on $\I^d$.

\bigskip

\begin{lemma}{} \label{A2} 
\begin{onelist}
\item 
The marginal distribution functions of $C_{(1,\aaa,\bbb)}$ and $C_{(2,\aaa,\bbb)}$ are identical.

\item 
$ C - 2p \, C_{(1,\aaa,\bbb)} + 2p \, C_{(2,\aaa,\bbb)} \in \CC $.
\end{onelist}
\end{lemma}{}

\bigskip

\begin{proof}
Assertion (1) is immediate from the definition and 
implies that the map $ C - 2p \, C_{(1,\aaa,\bbb)} + 2p \, C_{(2,\aaa,\bbb)} $ has uniform margins.
Now, we prove (2).
To this end, consider $ \uuu, \vvv \in \I^d $ with $\uuu \leq \vvv$.
Since $C$ is a copula and $ \aaa \leq \bbb $, we first obtain
$$
  \sum_{K \subseteq \{1,...,d\}} (-1)^{d-|K|} \, \big( C - 2p \, C_{(1,\aaa,\bbb)} \big) \big( \eeta_K(\uuu,\vvv) \big)
	  =  Q^C \big[ [\uuu,\vvv] \cap \I^d \backslash \big( [{\bf 0}, \aaa] \cup [\bbb, {\bf 1}] \big) \big]
	\geq 0
$$
and hence
\begin{eqnarray*}
  \lefteqn{\sum_{K \subseteq \{1,...,d\}} (-1)^{d-|K|} \, \big( C - 2p \, C_{(1,\aaa,\bbb)} + 2p \, C_{(2,\aaa,\bbb)} \big) \big( \eeta_K(\uuu,\vvv) \big)}
	\\
	& \geq & 
           2 p \, \sum_{K \subseteq \{1,...,d\}} (-1)^{d-|K|} \, \big( C_{(2,\aaa,\bbb)} \big) \big( \eeta_K(\uuu,\vvv) \big)
	\\*
	& \geq & 0
\end{eqnarray*}
where the last inequality follows from the fact that $C_{(2,\aaa,\bbb)}$ is a distribution function.
It remains to show that the identity
$ \big( C - 2p \, C_{(1,\aaa,\bbb)} + 2p \, C_{(2,\aaa,\bbb)} \big) (\eeta_k(\uuu,\zero)) = 0 $   
holds for every $ k\in\{1,...,d\} $ and every $ \uuu\in\I^d $, 
but this follows from the fact $C$ is a copula and $C_{(1,\aaa,\bbb)}$ and $C_{(2,\aaa,\bbb)}$ are distribution functions on $\I^d$.
\end{proof}

\bigskip

Motivated by Lemma \ref{A2}, 
we now put $ D :=  C - 2p \, C_{(1,\aaa,\bbb)} + 2p \, C_{(2,\aaa,\bbb)} $ and show that $ D \preceq C $ with $D \neq C$.

\bigskip

\begin{lemma}{} \label{A3}
\begin{onelist}
\item 
The inequalities
$ C_{(2,\aaa,\bbb)} (\uuu)
	\leq C_{(1,\aaa,\bbb)} (\uuu) $
and 
$ D(\uuu) \leq C(\uuu) $ hold for every $ \uuu \in \I^d $.

\item
There exists some $ \uuu \in \I^d $ satisfying
$ C_{(2,\aaa,\bbb)} (\uuu)
	< C_{(1,\aaa,\bbb)} (\uuu) $
and hence $ D(\uuu) < C(\uuu) $.

\item 
The inequalities
$ (\tau(C))_{(2,{\bf 1}-\bbb,{\bf 1}-\aaa)} (\uuu)
	\leq (\tau(C))_{(1,{\bf 1}-\bbb,{\bf 1}-\aaa)} (\uuu) $
and $ (\tau (D))(\uuu) \leq (\tau (C))(\uuu) $
hold for every $ \uuu \in \I^d $.

\item
We have $ D \preceq C $ with $ D \neq C $. 
\end{onelist}{}
\end{lemma}{}

\bigskip

\begin{proof}
To ease notation, for $\vvv \in \I^d$, we put $\vvv_1 := (v_2,\dots,v_d)$.
We first prove (1) and consider four cases.
\begin{hylist}
\item
First, assume that $u_1 < b_1$.
We then have 
\begin{eqnarray*}
  C_{(1,\aaa,\bbb)} (\uuu)
	&   =  & \frac{1}{2 \, p} \; C(\uuu \wedge \aaa) 
	\\
  C_{(2,\aaa,\bbb)} (\uuu)
	&   =  & \frac{1}{2 \, p^2} \; C(u_1 \wedge a_1, \aaa_1) \cdot  
					 Q^C \big[ [b_1,1] \times [\bbb_1,\uuu_1] \big]
\end{eqnarray*}
\begin{hylist}
\item 
If $u_i < b_i$ for some $i \in \{2,...,d\}$, 
then we obtain
$$
  C_{(2,\aaa,\bbb)} (\uuu)
	  =  \frac{1}{2 \, p^2} \; C(u_1 \wedge a_1, \aaa_1) \cdot  
					 Q^C \big[ [b_1,1] \times [\bbb_1,\uuu_1] \big]
		=  0
	\leq C_{(1,\aaa,\bbb)} (\uuu)
$$

\item
If $u_i \geq b_i \geq a_i$ for every $i \in \{2,...,d\}$, 
then we obtain
\begin{eqnarray*}
  C_{(2,\aaa,\bbb)} (\uuu)
	&   =  & \frac{1}{2 \, p^2} \; C(u_1 \wedge a_1, \aaa_1) \cdot  
					 Q^C \big[ [b_1,1] \times [\bbb_1,\uuu_1] \big]
	\\*
	& \leq & \frac{1}{2 \, p^2} \; C(u_1 \wedge a_1, \aaa_1) \cdot  Q^C \big[ [\bbb,{\bf 1}] \big]
	\\
	&   =  & \frac{1}{2 \, p} \; C(u_1 \wedge a_1, \aaa_1) 
	\\
	&   =  & C_{(1,\aaa,\bbb)} (\uuu)
\end{eqnarray*}
\end{hylist}

\item
Now, assume that $u_1 \geq b_1 \geq a_1$.
We then have 
\begin{eqnarray*}
  C_{(1,\aaa,\bbb)} (\uuu)
	&   =  & \frac{1}{2 \, p} \; \biggl( C (a_1, \uuu_1 \wedge \aaa_1) + Q^C \big[ [b_1, u_1] \times [\bbb_1,\uuu_1] \big] \biggr)
	\\
  C_{(2,\aaa,\bbb)} (\uuu)
  & = &	\frac{1}{2\,p^2} \; \biggl( C(\aaa) \cdot
				 Q^C \big[ [b_1,1] \times [\bbb_1,\uuu_1] \big]
				 + Q^C \big[ [b_1, u_1] \times [\bbb_1,{\bf 1}] \big] \cdot
				 C (a_1, \uuu_1 \wedge \aaa_1) \biggr)
\end{eqnarray*}
\begin{hylist}
\item
If $u_i < b_i$ for some $i \in \{2,...,d\}$, 
then we obtain
\begin{eqnarray*}
  C_{(2,\aaa,\bbb)} (\uuu)
  &   =  & \frac{1}{2\,p^2} \; Q^C \big[ [b_1, u_1] \times [\bbb_1,{\bf 1}] \big] \cdot
				 C (a_1, \uuu_1 \wedge \aaa_1) 
  \\
	& \leq & \frac{1}{2\,p^2} \; Q^C \big[ [\bbb,{\bf 1}] \big] \cdot
				 C (a_1, \uuu_1 \wedge \aaa_1)
	\\
	&   =  & \frac{1}{2\,p} \; C (a_1, \uuu_1 \wedge \aaa_1)
	\\
	&   =  & C_{(1,\aaa,\bbb)} (\uuu)
\end{eqnarray*}

\item
If $u_i \geq b_i \geq a_i$ for every $i \in \{2,...,d\}$, 
then we obtain
\begin{eqnarray*}
  C_{(2,\aaa,\bbb)} (\uuu)
	&   =  & \frac{1}{2 \, p^2} \; \biggl( C(\aaa) \cdot 
																 Q^C \big[ [b_1,1] \times [\bbb_1,\uuu_1] \big]
				   + Q^C \big[ [b_1,u_1] \times [\bbb_1,{\bf 1}] \big] \cdot C(\aaa) \biggr)
	\\
	&   =  & \frac{1}{2 \, p} \; \Big( Q^C \big[ [b_1,1] \times [\bbb_1,\uuu_1] \big]
					 + Q^C \big[ [b_1,u_1] \times [\bbb_1,{\bf 1}] \big]  \Big)
	\\
	&   =  & \frac{1}{2 \, p} \; \Big( Q^C \big[ [\bbb,\uuu] \big]
	         + Q^C \big[ [u_1,1] \times [\bbb_1,\uuu_1] \big]
					 + Q^C \big[ [b_1,u_1] \times [\bbb_1,{\bf 1}] \big] \Big)
	\\
	& \leq & \frac{1}{2 \, p} \; \Big( Q^C \big[ [\bbb,\uuu] \big]
	         + Q^C \big[ [u_1,1] \times [\bbb_1,{\bf 1}] \big]
					 + Q^C \big[ [b_1,u_1] \times [\bbb_1,{\bf 1}] \big] \Big)
	\\
	&   =  & \frac{1}{2 \, p} \; \Big( Q^C \big[ [\bbb,\uuu] \big]
	         + Q^C \big[ [\bbb,{\bf 1}] \big] \Big)
  \\
	&   =  & \frac{1}{2 \, p} \; \Big( p + Q^C \big[ [\bbb,\uuu] \big] \Big)
	\\
	&   =  & \frac{1}{2 \, p} \; \Big( C(\aaa) + Q^C \big[ [\bbb,\uuu] \big] \Big)
	\\
	&   =  & C_{(1,\aaa,\bbb)} (\uuu)
\end{eqnarray*}
\end{hylist}
\end{hylist}
This proves (1), and
(2) follows from the inequality
$$
  C_{(2,\aaa,\bbb)} (\aaa)
	= 0
	< \frac{1}{2}
	= C_{(1,\aaa,\bbb)} (\aaa)
$$
We now prove (3).
The inequality 
$ (\tau(C))_{(2,{\bf 1}-\bbb,{\bf 1}-\aaa)} (\uuu)
	\leq (\tau(C))_{(1,{\bf 1}-\bbb,{\bf 1}-\aaa)} (\uuu) $
for all $ \uuu \in \I^d $ is immediate from (1) and the fact that $ {\bf 1}-\bbb \leq {\bf 1}-\aaa $.
By \cite[Theorem 4.1]{fuc2014gamma},
we further have
\begin{eqnarray*}
  \big( \tau (D) \big) (\uuu)
	& = & \sum_{L \subseteq \{1,...,d\}} (-1)^{d-|L|} \, D \big( \eeta_L({\bf 1} - \uuu, {\bf 1}) \big)  
	\\*
	& = & \big( \tau (C) \big) (\uuu) +	
	      2p \, \sum_{L \subseteq \{1,...,d\}} (-1)^{d-|L|} \, 
				\Big( C_{(2,\aaa,\bbb)} \big( \eeta_L({\bf 1} - \uuu, {\bf 1}) \big) 
							- C_{(1,\aaa,\bbb)} \big( \eeta_L({\bf 1} - \uuu, {\bf 1}) \big) \Big)  
\end{eqnarray*}
for all $ \uuu \in \I^d $.
Moreover, applying \cite[Theorem 2.2]{fuc2016bic}, the identities
\begin{eqnarray*}
  \lefteqn{2p \, \sum_{L \subseteq \{1,...,d\}} (-1)^{d-|L|} \, C_{(1,\aaa,\bbb)} \big( \eeta_L({\bf 1} - \uuu, {\bf 1}) \big)}
	\\
  & = & Q^C \big[ [{\bf 1}-\uuu,{\bf 1}] \cap [{\bf 0}, \aaa] \big]
	      + Q^C \big[ [{\bf 1}-\uuu,{\bf 1}] \cap [\bbb, {\bf 1}] \big]
	\\
  & = & Q^{\tau(C)} \big[ [{\bf 0},\uuu] \cap [{\bf 1}-\aaa,{\bf 1}] \big]
	      + Q^{\tau(C)} \big[ [{\bf 0},\uuu] \cap [{\bf 0}, {\bf 1}-\bbb] \big]
	\\
	& = & 2p \, (\tau(C))_{(1,{\bf 1}-\bbb,{\bf 1}-\aaa)} (\uuu)
\end{eqnarray*}
and 
\begin{eqnarray*}
  \lefteqn{2p \, \sum_{L \subseteq \{1,...,d\}} (-1)^{d-|L|} \, C_{(2,\aaa,\bbb)} \big( \eeta_L({\bf 1} - \uuu, {\bf 1}) \big)}
	\\
  & = & \frac{1}{p} \;  Q^C \big[ [1-u_1,1] \times \I^{d-1} \cap [{\bf 0}, \aaa] \big] \cdot
				Q^C \big[ \I \times [{\bf 1}-\uuu_1,{\bf 1}] \cap [\bbb,{\bf 1}] \big]
  \\*
	&   & + \frac{1}{p} \; Q^C \big[ [1-u_1,1] \times \I^{d-1} \cap [\bbb,{\bf 1}] \big] \cdot
				Q^C \big[ \I \times [{\bf 1}-\uuu_1,{\bf 1}] \cap [{\bf 0}, \aaa] \big]
  \\
  & = & \frac{1}{p} \;  Q^{\tau(C)} \big[ [0,u_1] \times \I^{d-1} \cap [{\bf 1}-\aaa,{\bf 1}] \big] \cdot
				Q^{\tau(C)} \big[ \I \times [{\bf 0},\uuu_1] \cap [{\bf 0}, {\bf 1}-\bbb] \big] 
  \\*
	&   & + \frac{1}{p} \; Q^{\tau(C)}\big[ [0,u_1] \times \I^{d-1} \cap [{\bf 0}, {\bf 1}-\bbb] \big] \cdot
				Q^{\tau(C)} \big[ \I \times [{\bf 0},\uuu_1] \cap [{\bf 1}-\aaa,{\bf 1}] \big]
	\\*
	& = & 2p \, (\tau(C))_{(2,{\bf 1}-\bbb,{\bf 1}-\aaa)} (\uuu)
\end{eqnarray*}
hold for all $ \uuu \in \I^d $ and hence 
\begin{eqnarray*}
  \big( \tau (D) \big) (\uuu)
	&   =  & \big( \tau (C) \big) (\uuu) +	
	      2p \, \Big( (\tau(C))_{(2,{\bf 1}-\bbb,{\bf 1}-\aaa)} (\uuu)  - (\tau(C))_{(1,{\bf 1}-\bbb,{\bf 1}-\aaa)} (\uuu) \Big) 
	\\
	& \leq & \big( \tau (C) \big) (\uuu)
\end{eqnarray*}
for all $ \uuu \in \I^d $. This proves (3), and (4) is a consequence of (1), (2) and (3).
\end{proof}


\section*{Acknowledgements}
Jae Youn Ahn was supported by a National Research Foundation of Korea (NRF) grant funded by the Korean Government (NRF-2017R1D1A1B03032318).



{\footnotesize
\begin{thebibliography}{}

\bibitem[\protect\citeauthoryear{Ahn}{Ahn}{2015}]{Ahn11}
Ahn, J.~Y. (2015).
\newblock Negative dependence concept in copulas and the marginal free herd
  behavior index.
\newblock {\em J. Comput. Appl. Math.\/}~{\em 288}, 304--322.

\bibitem[\protect\citeauthoryear{Cap{\'e}ra{\`a} and Genest}{Cap{\'e}ra{\`a}
  and Genest}{1993}]{Caperaa}
Cap{\'e}ra{\`a}, P. and C.~Genest (1993).
\newblock Spearman's {$\rho$} is larger than {K}endall's {$\tau$} for
  positively dependent random variables.
\newblock {\em J. Nonparametr. Statist.\/}~{\em 2\/}(2), 183--194.

\bibitem[\protect\citeauthoryear{Cheung and Lo}{Cheung and Lo}{2014}]{Cheung4}
Cheung, K.~C. and A.~Lo (2014).
\newblock Characterizing mutual exclusivity as the strongest negative
  multivariate dependence structure.
\newblock {\em Insurance Math. Econom.\/}~{\em 55}, 180--190.

\bibitem[\protect\citeauthoryear{Dhaene and Denuit}{Dhaene and
  Denuit}{1999}]{Dhaene4}
Dhaene, J. and M.~Denuit (1999).
\newblock The safest dependence structure among risks.
\newblock {\em Insurance Math. Econom.\/}~{\em 25\/}(1), 11--21.

\bibitem[\protect\citeauthoryear{Dhaene, Denuit, Goovaerts, Kaas, and
  Vyncke}{Dhaene et~al.}{2002}]{Dhaene}
Dhaene, J., M.~Denuit, M.~J. Goovaerts, R.~Kaas, and D.~Vyncke (2002).
\newblock The concept of comonotonicity in actuarial science and finance: theory.
\newblock {\em Insurance Math. Econom.\/}~{\em 31\/}(1), 3--33.

\bibitem[\protect\citeauthoryear{Dhaene, Linders, Schoutens, and Vyncke}{Dhaene
  et~al.}{2014}]{daniel14}
Dhaene, J., D.~Linders, W.~Schoutens, and D.~Vyncke (2014).
\newblock A multivariate dependence measure for aggregating risks.
\newblock {\em J. Comput. Appl. Math.\/}~{\em 263},
  78--87.

\bibitem[\protect\citeauthoryear{Dolati and {\'U}beda-Flores}{Dolati and
  {\'U}beda-Flores}{2006}]{duf2006}
Dolati, A. and M.~{\'U}beda-Flores (2006).
\newblock On measures of multivariate concordance.
\newblock {\em J. Probab. Stat. Sci.\/}~{\em 4(2)}, 147--163.

\bibitem[\protect\citeauthoryear{Durante and Fern{\'a}ndez~S{\'a}nchez}{Durante
  and Fern{\'a}ndez~S{\'a}nchez}{2010}]{dufs2010}
Durante, F. and J.~Fern{\'a}ndez~S{\'a}nchez (2010).
\newblock Multivariate shuffles and approximation of copulas.
\newblock {\em Statist. Probab. Lett.\/}~{\em 80}, 1827--1834.

\bibitem[\protect\citeauthoryear{Durante and Sempi}{Durante and
  Sempi}{2016}]{dus2016}
Durante, F. and C.~Sempi (2016).
\newblock {\em Principles of Copula Theory}.
\newblock CRC Press, Boca Raton Fl.

\bibitem[\protect\citeauthoryear{Edwards, Mikusi{\'n}ski, and Taylor}{Edwards
  et~al.}{2004}]{emt2004}
Edwards, H.~H., P.~Mikusi{\'n}ski, and M.~D. Taylor (2004).
\newblock Measures of concordance determined by ${D}_{4}$--invariant copulas.
\newblock {\em Int. J. Math. Math. Sci.\/}~{\em 2004(70)}, 3867--3875.

\bibitem[\protect\citeauthoryear{Edwards, Mikusi{\'n}ski, and Taylor}{Edwards
  et~al.}{2005}]{emt2005}
Edwards, H.~H., P.~Mikusi{\'n}ski, and M.~D. Taylor (2005).
\newblock Measures of concordance determined by {$D\sb 4$}-invariant measures
  on {$(0,1)\sp 2$}.
\newblock {\em Proc. Amer. Math. Soc.\/}~{\em 133\/}(5), 1505--1513.

\bibitem[\protect\citeauthoryear{Fredricks and Nelsen}{Fredricks and
  Nelsen}{2007}]{Fredricks}
Fredricks, G.~A. and R.~B. Nelsen (2007).
\newblock On the relationship between {S}pearman's rho and {K}endall's tau for
  pairs of continuous random variables.
\newblock {\em J. Statist. Plann. Inference\/}~{\em 137\/}(7), 2143--2150.

\bibitem[\protect\citeauthoryear{Fuchs}{Fuchs}{2014}]{fuc2014gamma}
Fuchs, S. (2014).
\newblock Multivariate copulas: Transformations, symmetry, order and measures
  of concordance.
\newblock {\em Kybernetika\/}~{\em 50}, 725--743.

\bibitem[\protect\citeauthoryear{Fuchs}{Fuchs}{2016a}]{fuc2016bic}
Fuchs, S. (2016a).
\newblock A biconvex form for copulas.
\newblock {\em Depend. Model.\/}~{\em 4}, 63--75.

\bibitem[\protect\citeauthoryear{Fuchs}{Fuchs}{2016b}]{fuc2016moc}
Fuchs, S. (2016b).
\newblock Copula-induced measures of concordance.
\newblock {\em Depend. Model.\/}~{\em 4}, 205--214.

\bibitem[\protect\citeauthoryear{Fuchs, McCord, and Schmidt}{Fuchs
  et~al.}{2018}]{fumcsc2018}
Fuchs, S., Y.~McCord, and K.~D. Schmidt (2018).
\newblock Characterizations of copulas attaining the bounds of multivariate
  Kendall's tau.
\newblock {\em J. Optim. Theory Appl.\/}~{\em 178(2)}, 424--438.

\bibitem[\protect\citeauthoryear{Fuchs and Schmidt}{Fuchs and
  Schmidt}{2014}]{fus2014}
Fuchs, S. and K.~D. Schmidt (2014).
\newblock Bivariate copulas: Transformations, asymmetry and measures of
  concordance.
\newblock {\em Kybernetika\/}~{\em 50}, 109--125.

\bibitem[\protect\citeauthoryear{Genest, Ne{\v s}lehov{\'a}, and
  Ben~Ghorbal}{Genest et~al.}{2011}]{genebg2011}
Genest, C., J.~Ne{\v s}lehov{\'a}, and N.~Ben~Ghorbal (2011).
\newblock Estimators based on Kendall's tau in multivariate copula models.
\newblock {\em Aust. N. Z. J. Stat.\/}~{\em 53}, 157--177.

\bibitem[\protect\citeauthoryear{H{\"u}rlimann}{H{\"u}rlimann}{2003}]{Hurlimann}
H{\"u}rlimann, W. (2003).
\newblock Hutchinson-{L}ai's conjecture for bivariate extreme value copulas.
\newblock {\em Statist. Probab. Lett.\/}~{\em 61\/}(2), 191--198.

\bibitem[\protect\citeauthoryear{Joe}{Joe}{1990}]{joe1990}
Joe, H. (1990).
\newblock Multivariate concordance.
\newblock {\em J. Multivariate Anal.\/}~{\em 35}, 12--30.

\bibitem[\protect\citeauthoryear{Lee and Ahn}{Lee and Ahn}{2014}]{Ahn7}
Lee, W. and J.~Y. Ahn (2014).
\newblock On the multidimensional extension of countermonotonicity and its
  applications.
\newblock {\em Insurance Math. Econom.\/}~{\em 56}, 68--79.

\bibitem[\protect\citeauthoryear{Lee, Cheung, and Ahn}{Lee et~al.}{2017}]{Ahn9}
Lee, W., K.~C. Cheung, and J.~Y. Ahn (2017).
\newblock Multivariate countermonotonicity and the minimal copulas.
\newblock {\em J. Comput. Appl. Math.\/}~{\em 317}, 589--602.

\bibitem[\protect\citeauthoryear{Lux and Papapantoleon}{Lux and
  Papapantoleon}{2017}]{luxpa2017}
Lux, T. and A.~Papapantoleon (2017).
\newblock Improved {F}r\'echet-Hoeffding bounds on $d-$copulas and applications
  in model-free finance.
\newblock {\em Ann. Appl. Probab.\/}~{\em 27\/}(6), 3633--3671.

\bibitem[\protect\citeauthoryear{M{\"u}ller and Scarsini}{M{\"u}ller and
  Scarsini}{2000}]{musc2000}
M{\"u}ller, A. and M.~Scarsini (2000).
\newblock Some remarks on the supermodular order.
\newblock {\em J. Multivariate Anal.\/}~{\em 73}, 107--119.

\bibitem[\protect\citeauthoryear{Nelsen}{Nelsen}{2002}]{nel2002}
Nelsen, R.~B. (2002).
\newblock Concordance and copulas: A survey.
\newblock In C.~M. Cuadras, J.~Fortiana, and J.~A. Rodriguez-Lallena (Eds.),
  {\em Distributions with Given Marginals and Statistical Modelling}, pp.\
  169--177. 
Kluwer Academic Publishers, Dordrecht.

\bibitem[\protect\citeauthoryear{Nelsen}{Nelsen}{2006}]{nel2006}
Nelsen, R.~B. (2006).
\newblock {\em An Introduction to Copulas.\/} Second edition.
\newblock Springer, New York.

\bibitem[\protect\citeauthoryear{Nelsen and {\'U}beda-Flores}{Nelsen and
  {\'U}beda-Flores}{2012}]{neuf2012}
Nelsen, R.~B. and M.~{\'U}beda-Flores (2012).
\newblock Directional dependence in multivariate distributions.
\newblock {\em Ann. Inst. Stat. Math.\/}~{\em 64}, 677--685.

\bibitem[\protect\citeauthoryear{Preischl}{Preischl}{2016}]{prei2016}
Preischl, M. (2016).
\newblock Bounds on integrals with respect to multivariate copulas..
\newblock {\em Depend. Model.\/}~{\em 4}, 277--287.

\bibitem[\protect\citeauthoryear{Puccetti and Wang}{Puccetti and
  Wang}{2015}]{ruodu3}
Puccetti, G. and R.~Wang (2015).
\newblock Extremal dependence concepts.
\newblock {\em Statist. Sci.\/}~{\em 30}, 485--517.

\bibitem[\protect\citeauthoryear{Rudin}{Rudin}{1976}]{Rudin1976}
Rudin, W. (1976).
\newblock {\em Principles of Mathematical Analysis.\/} Third edition.
\newblock Chapman \& Hall, London.

\bibitem[\protect\citeauthoryear{R{\"u}schendorf and Uckelmann}{R{\"u}schendorf
  and Uckelmann}{2002}]{Ruschendorf2}
R{\"u}schendorf, L. and L.~Uckelmann (2002).
\newblock Variance minimization and random variables with constant sum.
\newblock In {\em Distributions with given marginals and statistical
  modelling}, pp.\  211--222. 
Kluwer Academic Publishers, Dordrecht. 

\bibitem[\protect\citeauthoryear{Scarsini}{Scarsini}{1984}]{sca1984}
Scarsini, M. (1984).
\newblock On measures of concordance.
\newblock {\em Stochastica\/}~{\em 8(3)}, 201--218.

\bibitem[\protect\citeauthoryear{Schreyer, Paulin, and Trutschnig}{Schreyer
  et~al.}{2017}]{schpautru2017}
Schreyer, M., R.~Paulin, and W.~Trutschnig (2017).
\newblock On the exact region determined by Kendall's $\tau$ and Spearman's
  $\rho$.
\newblock {\em J. R. Stat. Soc. Ser. B Stat. Methodol.\/}~{\em 79}, 613--633.

\bibitem[\protect\citeauthoryear{Schweizer and Sklar}{Schweizer and
  Sklar}{1983}]{Schweizer}
Schweizer, B. and A.~Sklar (1983).
\newblock {\em Probabilistic metric spaces}.
\newblock North-Holland, New York.

\bibitem[\protect\citeauthoryear{Tankov}{Tankov}{2011}]{Tankov}
Tankov, P. (2011).
\newblock Improved {F}r\'echet bounds and model-free pricing of multi-asset
  options.
\newblock {\em J. Appl. Probab.\/}~{\em 48\/}(2), 389--403.

\bibitem[\protect\citeauthoryear{Taylor}{Taylor}{2007}]{tay2007}
Taylor, M.~D. (2007).
\newblock Multivariate measures of concordance.
\newblock {\em Ann. Inst. Statist. Math.\/}~{\em 59}, 789--806.

\bibitem[\protect\citeauthoryear{Taylor}{Taylor}{2016}]{tay2016}
Taylor, M.~D. (2016).
\newblock Multivariate measures of concordance for copulas and their marginals.
\newblock {\em Depend. Model.\/}~{\em 4\/}, 224--236.

\bibitem[\protect\citeauthoryear{{\'U}beda-Flores}{{\'U}beda-Flores}{2005}]{ubf2005}
{\'U}beda-Flores, M. (2005).
\newblock Multivariate versions of Blomqvist's beta and Spearman's footrule.
\newblock {\em Ann. Inst. Statist. Math.\/}~{\em 57}, 781--788.

\bibitem[\protect\citeauthoryear{Wang and Wang}{Wang and Wang}{2011}]{wawa2011}
Wang, B. and R.~Wang (2011).
\newblock The complete mixability and convex minimization problems with
  monotone marginal densities.
\newblock {\em J. Multivariate Anal.\/}~{\em 102 (10)}, 1344--1360.

\bibitem[\protect\citeauthoryear{Wang and Wang}{Wang and Wang}{2016}]{Ruodu4}
Wang, B. and R.~Wang (2016).
\newblock Joint mixability.
\newblock {\em Math. Oper. Res.\/}~{\em 41\/}(3), 808--826.

\end{thebibliography}
}
\end{document}